\documentclass{owrart}

\newtheorem{thm}{Theorem}[section]

\newtheorem{ques}[thm]{Question}

\newcommand{\be}{\begin{equation}}

\newcommand{\ee}{\end{equation}}

\newcommand{\N}{\mathbb{N}}

\newcommand{\R}{\mathbb{R}}

\DeclareMathOperator{\Vol}{Vol}
\DeclareMathOperator{\Diam}{Diam}
\DeclareMathOperator{\Area}{Area}

\newcommand{\Sp}{\mathbb{S}}      
\newcommand{\Tor}{\mathbb{T}} 

 \usepackage{graphicx}
 \usepackage{tikz-cd}
\usepackage{tikz}

\begin{document}


\begin{talk}[E. Bryden, D. Kazaras, R. Perales, and C. Sormani]{Brian Allen}
{Scalar Curvature Stability: Tools, Theorems, and Questions}

\noindent

\vspace{-0.25cm}

\section{Introduction}
Although scalar curvature is the simplest curvature invariant, our understanding of scalar curvature has not matured to the same level as Ricci or sectional curvature. Despite this fact, many rigidity phenomenon have been established which give some of the strongest insights into scalar curvature. Important examples include Geroch's conjecture, the positive mass theorem, and Llarull's theorem. In order to further understand scalar curvature we ask corresponding geometric stability questions, where the hypotheses of the rigidity phenomenon are relaxed, and one would like to show that Riemannian manifolds which satisfy the relaxed conditions are close to the rigid objects in some topology. In this talk we will survey what is known for scalar curvature stability, discuss what the questions are in this area, and introduce important tools which have been useful so far. 

\vspace{-0.25cm}

\section{Scalar Curvature Phenomenon}

When exploring sequences of Riemannian manifolds with scalar curvature lower bounds we see that splines, bubbles, and drawstrings, depicted in Figure \ref{fig-phenomenon}, are persistent phenomenon which need to either be ruled out by making an assumption or addressed by the choice of topology one makes when choosing an appropriate notion of convergence. Examples with splines and bubbles are constructed using Gromov-Lawson tunnels \cite{Gromov-Lawson-torus} with estimates. Careful constructions of this type have been carried out by Basilio, Dodziuk, and Sormani \cite{BDS}, Basilio and Sormani \cite{Basilio-Sormani},  Dodziuk \cite{Dodziuk}, Basilio, Kazaras and Sormani \cite{Basilio-Kazaras-Sormani}, and Sweeney \cite{Sweeney}. If one allows a sequence of manifolds which are non-diffeomorphic to the limit manifold then Basilio and Sormani \cite{Basilio-Sormani} have also shown that sewing is possible, something we will mostly avoid in this note by fixing the topology. The existence of drawstrings in dimensions $n \ge 4$ with scalar curvature bounds was constructed by Lee, Naber, and Neumayer \cite{LNN} and in dimension $n=3$ by Kazaras and Xu \cite{KK}.

\begin{figure}[h] \label{fig-phenomenon}
   \center{\includegraphics[width=.6\textwidth]{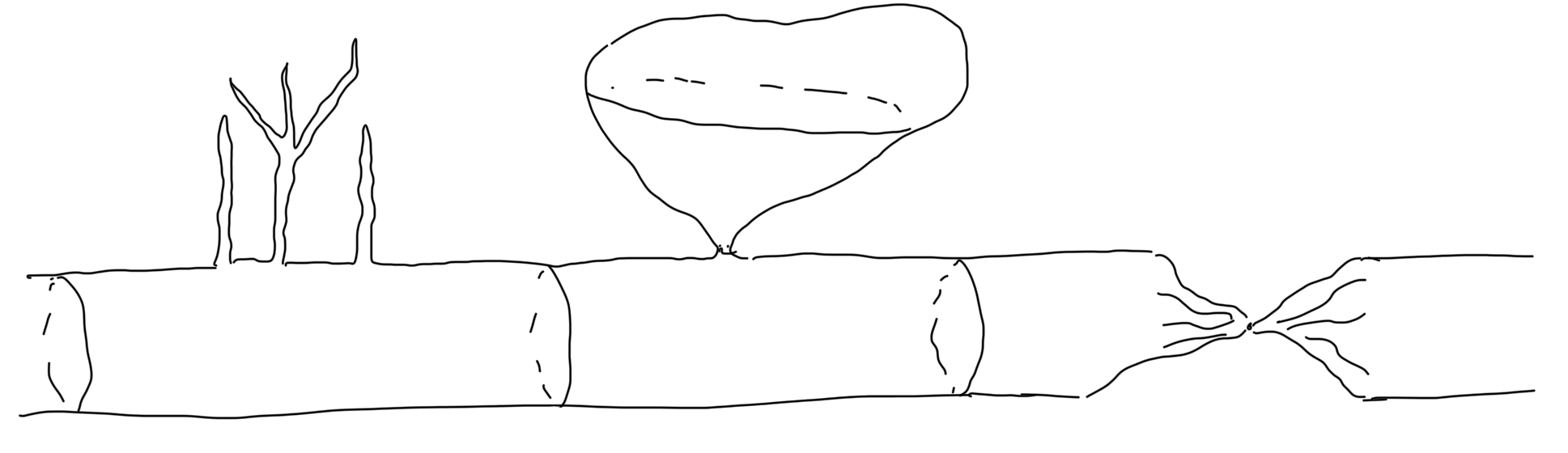}}
   \caption{From left to right, splines, bubbling, and a drawstring where a circle has been pulled to almost a point. Along such a sequence the splines would become arbitrarily thin, the neck of the bubble would pinch to a point, and the circle would be pulled to a point by the drawstring.}
\end{figure}

Due to the presence of splines in all scalar curvature stability problems we know that Gromov-Hausdorff stability is not appropriate for these problems since a sequence of manifolds with increasingly many splines cannot converge in the Gromov-Hausdorff sense by Gromov's compactness theorem \cite{BBI}. Hence we are left to look for other notions of convergence for sequences of Riemannian manifolds which are resilient to the phenomenon depicted in Figure \ref{fig-phenomenon}. In this note we will give examples of three different choices of convergence one can use: $d_p$ convergence (see Lee, Naber, Neumayer \cite{LNN} for the definition), Gromov-Hausdorff convergence of the restricted length metric to a good set, or Sormani-Wenger Intrinsic Flat (SWIF) convergence (see Sormani and Wenger \cite{SW-JDG} for the definition). We emphasize that there is no clear hierarchy between the three choices of convergence of Riemannian manifolds we will discuss in this note since each has drawbacks in how they handle the scalar curvature phenomenon depicted in Figure \ref{fig-phenomenon}. Hence it is the philosophy of the author that our goal should be to prove all three notions of convergence, under different assumptions, for most scalar curvature rigidity results. For example, one would rule out splines and bubbles if one wants to show $d_p$ stability and one would rule out bubbles and drawstrings if one would like to show SWIF stability. In concert, the three stability results would tell us the most about sequences of Riemannian manifolds with scalar curvature bounds, as will be discussed throughout the rest of this note.

\vspace{-0.15cm}

\section{Geroch Conjecture Stability}

\begin{thm}[Schoen and Yau \cite{Schoen-Yau-min-surf}, Gromov and Lawson \cite{Gromov-Lawson-torus}]
If $(\Tor^n,g)$ is a Riemannian metric with non-negative scalar curvature then $(\Tor^n,g)$ is isometric to a flat torus.
\end{thm}

Various special cases of Geroch stability have already been addressed. Single and double warped products were addressed by Allen, Vazquez, Parise, Payne, and Wang \cite{AHMPPW}, the graph case was studied by Pacheco, Ketterer, and Perales \cite{PKP19}, and conformal cases by Allen \cite{Allen-Conformal-Torus}, and Chu and Lee \cite{Chu-Man-Chun22}, and under isoperimetric bounds and integral Ricci bounds by Allen, Bryden, and Kazaras \cite{Allen-Bryden-Kazaras-PMT-Stability}. In these special cases, many different notions of convergence were considered. In general, since splines, bubbles, and drawstrings are possible in the case of Geroch stability, one needs to add a hypothesis to remove two of the three scalar curvature phenomenon to show $d_p$ or SWIF stability or one needs to cut these phenomenon out altogether. In the following stability result, the authors choose to assume an entropy bound, which rules out splines and bubbles, and show $d_p$ convergence to a flat torus. One can find many examples in \cite{LNN} where it is demonstrated that $d_p$ convergence is resilient to drawstrings. See \cite{Allen-Bryden-d_p} for an explanation of the fact that $d_p$ convergence is not well suited in the presence of splines and bubbles.



\begin{thm}[Lee, Naber, and Neumayer \cite{LNN}]\label{thm-LNN}
 Fix $n \ge 2$ and $p \ge n+1$. There exists a $\delta=\delta(n,p)>0$ and $V_0=v_0(n,p)>0$ such that the following holds: For any $V>V_0$ and $(\mathbb{T}^n,g_j)$  a sequence of Riemannian tori such that
   \begin{align}
     R_{g_j} \ge -\frac{1}{j}, \qquad \Vol(\mathbb{T}^n,g_j) \le V_0, \qquad \nu(\mathbb{T}^n,g_j)\ge -\delta,
   \end{align}
   a subsequence of $(\mathbb{T}^n,g_j)$ converges to a flat torus $(\mathbb{T}^n,g_F)$ in the $d_p$ sense.

\end{thm}

\begin{ques}
Given that the entropy bound of Theorem \ref{thm-LNN} removes the possibility of splines and bubbles forming along a sequence, can we formulate and prove Geroch stability where one removes the possibility of bubbles and drawstrings and shows volume preserving SWIF convergence to a flat torus?
\end{ques}
\vspace{-0.25cm}

\section{Positive Mass Theorem Stability}

\begin{thm}[Schoen and Yau \cite{SchoenYauPMT}, Witten \cite{Witten-PMT}]

If $(M^3,g)$ is an asympotically flat manifold with non-negative scalar curvature then the ADM mass is non-negative. If the ADM mass is zero then $(M^3,g)$ is isometric to Euclidean space.
\end{thm}

The stability of the positive mass theorem has been studied in many special cases. Assuming the existence of a smooth IMCF by Allen \cite{Allen18,Allen-AsymHypPMT,Allen-SobolevPMT, Allen-SWIFPMT}, in terms of the Brown-York Mass for graphs in Euclidean space by Alaee, Cabrera Pacheco, and McCormick \cite{Alaee-Cabrera-McCormick}, metrics conformal to Euclidean space by Corvino \cite{Corvino-PMT}, Lee \cite{Lee-PMT-Conformal}, axially symmetric metrics by Bryden \cite{Bryden20}, using spinors by Finster and Bray \cite{Finster-Bray-99}, Finster and Kath \cite{Finster-Kath-02}, and Finster \cite{Finster-09},  in the geometrostatic setting by Stavrov-Allen and Sormani \cite{Stavrov-Sormani}, in the rotationally symmetric setting by Lee and Sormani \cite{LeeSormani1}, assuming a lower Ricci bound by Kazaras, Lee, and Khuri \cite{Lee-Kazaras-Khuri-PMT-Stability}, and under isoperimetric bounds and integral Ricci bounds by Allen, Bryden, and Kazaras \cite{Allen-Bryden-Kazaras-PMT-Stability}. In these special cases, many different notions of convergence were considered.

In the following result, the authors address a conjecture by Huisken  and Ilmanen \cite{Huisken-Ilmanen} where the idea is to remove all of the scalar curvature phenomenon depicted in Figure \ref{fig-phenomenon} and show that the remaining set converges to Euclidean space.

\begin{thm}[Dong and Song \cite{Dong-Song}]\label{thm-Dong-Song}
Let $(M_i^3,g_i)$ be a sequence of asymptotically flat Riemannian manifolds with non-negative scalar curvature and suppose that the ADM mass $m_{ADM}(M_i,g_i) \searrow 0$. Then for all $i \in \N$ and each end in $M_i$, there is a domain $Z_i \subset M_i$ with smooth boundary so that $|\partial Z_i|_{g_i} \searrow 0$, $M_i \setminus Z_i$ contains the given end, and
\begin{align}
    (M_i\setminus Z_i,\hat{d}_{g_i}, p_i) \rightarrow (\R^3, d_{\mathbb{E}^3},0),
\end{align}
in the pointed measured Gromov-Hausdorff sense, where $p_i \in M_i\setminus Z_i$ is any choice of base point, and $\hat{d}_{g_i}$ is the restricted length metric on $M_i \setminus Z_i$ induced by $g_i$.
\end{thm}

\begin{ques}
   Can one formulate and prove a version of Geroch stability in the style of Theorem \ref{thm-Dong-Song}?
\end{ques}

\begin{ques}
   In light of Theorem \ref{thm-Dong-Song}, one can ask to say more about what is happening on the bad set $Z_i$. For instance, is $Z_i$ made up of splines, bubbles, drawstrings, and sewing only? One way of providing some evidence in this direction is to formulate a stability conjecture where one assumes a condition which rules out splines and bubbles and concludes  $d_p$ stability and another version of a stability conjecture where one rules out bubbles, drawstrings, and sewing which concludes volume preserving SWIF convergence. In concert, these three stability results would be giving strong evidence that the bad sets $Z_i$ are made up of only the known phenomenon, reinforcing the belief that these are the only phenomenon which one needs to worry about when discussing sequences with scalar curvature lower bounds. 
\end{ques}

\vspace{-0.35cm}

\section{Llarull Stability}

Here we state a slightly less general version of Llarull's rigidity theorem.

\begin{thm}[Llarull \cite{Llarull}]\label{thm-Llarull}

If $(\Sp^n,g)$ is a Riemannian metric so that $g \ge g_{\Sp^n}$, where $g_{\Sp^n}$ is the round sphere, whose scalar curvature $R_g \ge n(n-1)$ then $(\Sp^n,g)$ is isometric to $(\Sp^n,g_{\Sp^n})$.
\end{thm}

In this case, when we consider the stability of Theorem \ref{thm-Llarull} we notice that drawstrings are not possible due to the metric lower bound assumption which is necessary in the rigidity theorem. Hence, volume preserving SWIF convergence is most appropriate for stability in this case. In fact, any time that rigidity is phrased in terms of a $1-$Lipschitz map $F:(M,g)\rightarrow (N,h)$ we claim that volume preserving Sormani-Wenger Intrinsic Flat convergence will be the correct choice for proving the corresponding stability result.

\begin{thm}[Allen, Bryden, and Kazaras \cite{ABKLLarull}]\label{thm-LlarullStability}
        Let $V,D,\overline{m},\Lambda>0$. If a sequence $(\Sp^3,g_i)$ of Riemannian $3$-spheres satisfies
        \begin{align}
      g_i&\geq g_{\Sp^3},
            \qquad \Vol(\Sp^3,g_i)\leq V,
            \\ \Diam(\Sp^3,g_i)&\le D,
            \qquad \inf_{\Omega\subset\Sp^3}\tfrac{\Area(\partial\Omega,g_i)}{\min(\Vol(\Omega,g_i),\Vol(\Sp^3\setminus\Omega,g_i))}\geq\Lambda,
        \end{align}
        and
        \begin{align}
\left\|\left(6-R_{g_i}\right)^{+}\right\|_{L^{2}(g_i)}^{1/2} \rightarrow 0,
        \end{align}
        then it converges in the volume preserving SWIF sense:
        \begin{align}
            d_{\mathcal{VF}}((\Sp^3,g_i),(\Sp^3,g_{\mathbb{S}^3}))\rightarrow 0.
        \end{align}
    \end{thm}

    The proof uses two main tools from the literature: spacetime harmonic functions developed by Hirsch, Kazaras, and Khuri \cite{hirsch-kazaras-khuri} and used to give a new proof of Llarull's theorem by Hirsch, Kazaras, Khuri, and  Zang \cite{Hirsch-Kazaras-Khuri-Zhang}, and the Volume Above Distance Below (VADB) theorem of Allen,  Perales, and Sormani \cite{Allen-Perales-Sormani-VADB}. Shortly after Theorem \ref{thm-LlarullStability} was established, a more general result was established using the spinor formulas which Llarull used to establish rigidity and the VADB theorem.

    \begin{thm}[Hirsch and Zhang \cite{HZ}]\label{thm-HZ}
Let $V,D,\overline{m},\Lambda>0$. If a sequence $(\Sp^n,g_i)$, $n \ge 3$ of Riemannian $n$-spheres satisfies
        \begin{align}
      g_i\geq g_{\Sp^n},
            \qquad  \Diam(\Sp^n,g_i)&\le D, \qquad\inf_{u \in W^{1,2}(\Sp^n)}\tfrac{\|\nabla u\|_{L^2(\Sp^n)}^2}{\inf_{k \in \R}\|u-k\|_{L^2(\Sp^n)}^2}\geq\Lambda,
        \end{align}
        and
        \begin{align}
R_{g_i} \ge n(n-1)-\frac{1}{i}
        \end{align}
        then it converges in the volume preserving SWIF sense:
        \begin{align}
            d_{\mathcal{VF}}((\Sp^n,g_i),(\Sp^n,g_{\mathbb{S}^n}))\rightarrow 0.
        \end{align}

\end{thm}

It should be noted that Hirsch and Zhang prove two other versions of Llarull stability, one with an $L^p$ assumption on the scalar curvature below $n(n-1)$ and another result with a good set bad set decomposition similar to Theorem \ref{thm-Dong-Song}.

\begin{ques}
Can one use spinors in order to establish Geroch stability and positive mass theorem stability in dimensions $n \ge 3$?
\end{ques}

\begin{ques}
Can one generalize Theorem \ref{thm-HZ} to the case of the results of Goette and Semmelmann \cite{Goette-Semmelmann}?
\end{ques}
\vspace{-0.25cm}

\section{Volume Above Distance Below Theorem}

We end this note by stating and discussing the VADB theorem which has proved useful so far in proving stability of Llarull's theorem as well as many special cases of Geroch stability and positive mass theorem stability.

\begin{thm}[Allen, Perales, and Sormani \cite{Allen-Perales-Sormani-VADB}]\label{thm-VADB} 
Suppose we have a fixed  closed, oriented, Riemannian manifold, $(M^n,g_0)$, and a sequence of  continuous Riemannian manifolds $(M,g_j)$ such that
\be\label{eq-thmbetterd}
g_j(v,v) \ge \left( 1-C_j\right)g_0(v,v), \quad \forall p \in M,  v \in T_pM, \quad C_j \searrow 0,
\ee
and a uniform upper bound on diameter, $\Diam(M_j) \le D_0$, and volume convergence
\be
\Vol(M,g_j) \to \Vol(M,g_0),
\ee
then we find volume preserving SWIF convergence $d_{\mathcal{VF}}((M,g_j), (M,g_0)) \rightarrow 0$.

\end{thm}

It should also be noted that Allen and Perales \cite{Allen-Perales} have proved a version of Theorem \ref{thm-VADB} for manifolds with boundary. Hence one can use the VADB theorem with boundary to discuss stability of rigidity results on manifolds with boundary.

\begin{ques}
 Can one prove a version of Theorem \ref{thm-VADB} for sequences which are not diffeomorphic to the limit manifold? Addressing this question would allow us to extend Theorem \ref{thm-LlarullStability} and Theorem \ref{thm-HZ} to $1-$Lipschitz maps from manifolds which are not diffeomorphic to the sphere.
\end{ques}

\begin{ques}
    Theorem \ref{thm-VADB} motivates a notion of convergence for sequences of Riemannian manifolds. We say that a sequence of Riemannian manifolds $(M,g_j)$ converges to a Riemannian manifold $(M,g_{\infty})$ in the VADB sense if $\Diam(M,g_j) \le D$, $\Vol(M,g_j) \rightarrow \Vol(M,g_{\infty})$ and $g_j(v,v) \ge (1-C_j)g_{\infty}(v,v)$, $\forall p \in M,  v \in T_pM, C_j \searrow 0$. If we consider a sequence of Riemannian manifolds $(M,g_j)$ so that $R_{g_j} \ge 0$, which converge in the VADB sense to $(M,g_{\infty})$, then is it the case that $R_{g_{\infty}}\ge 0$? It was shown by Gromov \cite{GroD} and Bamler \cite{Bamler} that if $g_j \rightarrow g_{\infty}$ in the $C^0$ sense then the answer to the question is yes. It was shown by Lee and Topping \cite{Lee-Topping} that if $d_{g_j} \rightarrow d_{g_{\infty}}$ in the uniform sense then the answer is no. Considering this question for VADB convergence is interesting because VADB convergence lies in between $C^0$ convergence of Riemannian metrics and uniform convergence of distance functions. It is also interesting since VADB convergence allows for the presence of splines along the sequences which persist in the limit.  \end{ques}


  \bibliographystyle{plain}
    \bibliography{bibliography}

\begin{thebibliography}{10}

\bibitem{Alaee-Cabrera-McCormick}
Aghil Alaee, Armando~J. Cabrera~Pacheco, and Stephen McCormick.
\newblock Stability of a quasi-local positive mass theorem for graphical
  hypersurfaces of {E}uclidean space.
\newblock {\em Trans. Amer. Math. Soc.}, 374(5):3535--3555, 2021.

\bibitem{Allen-Conformal-Torus}
B.~Allen.
\newblock Almost non-negative scalar curvature on riemannian manifolds
  conformal to tori.
\newblock {\em Journal of Geometric Analysis}, 31:11190--11213, 2021.

\bibitem{AHMPPW}
B.~Allen, L.~Hernandez-Vazquez, D.~Parise, A.~Payne, and S.~Wang.
\newblock Warped tori with almost non-negative scalar curvature.
\newblock {\em Geometriae Dedicata}, 2018.

\bibitem{Allen18}
Brian Allen.
\newblock Imcf and the stability of the pmt and rpi under $l^2$ convergence.
\newblock {\em Annales Henri Poincar\'e}, 19(1), 2017.

\bibitem{Allen-AsymHypPMT}
Brian Allen.
\newblock Stability of the {PMT} and {RPI} for asymptotically hyperbolic
  manifolds foliated by {IMCF}.
\newblock {\em J. Math. Phys.}, 59(8):082501, 18, 2018.

\bibitem{Allen-SobolevPMT}
Brian Allen.
\newblock Sobolev stability of the positive mass theorem and {R}iemannian
  {P}enrose inequality using inverse mean curvature flow.
\newblock {\em Gen. Relativity Gravitation}, 51(5):Paper No. 59, 32, 2019.

\bibitem{Allen-SWIFPMT}
Brian Allen.
\newblock Inverse mean curvature flow and the stability of the positive mass
  theorem.
\newblock {\em Communications in Analysis and Geometry}, 31, 2023.

\bibitem{Allen-Bryden-d_p}
Brian Allen and Edward Bryden.
\newblock Exploring a modification of $d_p$ convergence, 2023.
\newblock arXiv:2311.13450.

\bibitem{Allen-Bryden-Kazaras-PMT-Stability}
Brian Allen, Edward Bryden, and Demetre Kazaras.
\newblock Stability of the positive mass theorem and torus rigidity theorems
  under integral curvature bounds, 2022.
\newblock arXiv:2210.04340.

\bibitem{ABKLLarull}
Brian Allen, Edward Bryden, and Demetre Kazaras.
\newblock On the stability of llarull's theorem in dimension three, 2023.
\newblock arXiv:2305.18567.

\bibitem{Allen-Perales}
Brian Allen and Raquel Perales.
\newblock Intrinsic flat stability of manifolds with boundary where volume
  converges and distance is bounded below, 2021.
\newblock arXiv:2006.13030.

\bibitem{Allen-Perales-Sormani-VADB}
Brian Allen, Raquel Perales, and Christina Sormani.
\newblock Volume above distance below.
\newblock {\em To appear in the Journal of Differential Geometry}.

\bibitem{Bamler}
Richard~H. Bamler.
\newblock A {R}icci flow proof of a result by {G}romov on lower bounds for
  scalar curvature.
\newblock {\em Math. Res. Lett.}, 23(2):325--337, 2016.

\bibitem{BDS}
J.~Basilio, J.~Dodziuk, and C.~Sormani.
\newblock Sewing {R}iemannian manifolds with positive scalar curvature.
\newblock {\em J. Geom. Anal.}, 28(4):3553--3602, 2018.

\bibitem{Basilio-Kazaras-Sormani}
J.~Basilio, D.~Kazaras, and C.~Sormani.
\newblock An intrinsic flat limit of {R}iemannian manifolds with no geodesics.
\newblock {\em Geom. Dedicata}, 204:265--284, 2020.

\bibitem{Basilio-Sormani}
J.~Basilio and C.~Sormani.
\newblock Sequences of three dimensional manifolds with positive scalar
  curvature.
\newblock 2019.

\bibitem{Finster-Bray-99}
Hubert Bray and Felix Finster.
\newblock Curvature estimates and the positive mass theorem.
\newblock {\em Comm. Anal. Geom}, pages 291--306, 2002.

\bibitem{Bryden20}
Edward Bryden.
\newblock Stability of the positive mass theorem for axisymmetric manifolds.
\newblock {\em Pacific journal of mathematics}, 305(1):89--152, 2020.

\bibitem{BBI}
D.~Burago, Y.~Burago, and S.~Ivanov.
\newblock {\em A course in metric geometry}, volume~33 of {\em Graduate Studies
  in Mathematics}.
\newblock American Mathematical Society, Providence, RI.

\bibitem{Chu-Man-Chun22}
Jianchun Chu and Man-Chun Lee.
\newblock K{\"a}hler tori with almost non-negative scalar curvature.
\newblock {\em Communications in Contemporary Mathematics}, 0(0):2250030, 2022.

\bibitem{Corvino-PMT}
Justin Corvino.
\newblock A note on asymptotically flat metrics on {${\Bbb R}^3$} which are
  scalar-flat and admit minimal spheres.
\newblock {\em Proc. Amer. Math. Soc.}, 133(12):3669--3678, 2005.

\bibitem{Dodziuk}
J\'{o}zef Dodziuk.
\newblock Gromovl-{L}awson tunnels with estimates.
\newblock In {\em Analysis and geometry on graphs and manifolds}, volume 461 of
  {\em London Math. Soc. Lecture Note Ser.}, pages 55--65. Cambridge Univ.
  Press, Cambridge, 2020.

\bibitem{Dong-Song}
Conghan Dong and Antoine Song.
\newblock Stability of euclidean 3-space for the positive mass theorem, 2023.
\newblock arXiv:2302.07414.

\bibitem{Finster-09}
Felix Finster.
\newblock A level set analysis of the {W}itten spinor with applications to
  curvature estimates.
\newblock {\em Math. Res. Lett.}, 16(1):41--55, 2009.

\bibitem{Finster-Kath-02}
Felix Finster and Ines Kath.
\newblock Curvature estimates in asymptotically flat manifolds of positive
  scalar curvature.
\newblock {\em Comm. Anal. Geom.}, 10(5):1017--1031, 2002.

\bibitem{Goette-Semmelmann}
S.~Goette and U.~Semmelmann.
\newblock Scalar curvature estimates for compact symmetric spaces.
\newblock {\em Differential Geom. Appl.}, 16(1):65--78, 2002.

\bibitem{Gromov-Lawson-torus}
Mikhael Gromov and H.~Blaine Lawson, Jr.
\newblock Spin and scalar curvature in the presence of a fundamental group.
  {I}.
\newblock {\em Ann. of Math. (2)}, 111(2):209--230, 1980.

\bibitem{GroD}
Misha Gromov.
\newblock Dirac and {P}lateau billiards in domains with corners.
\newblock {\em Cent. Eur. J. Math.}, 12(8):1109--1156, 2014.

\bibitem{hirsch-kazaras-khuri}
Sven Hirsch, Demetre Kazaras, and Marcus Khuri.
\newblock Spacetime harmonic functions and the mass of 3-dimensional
  asymptotically flat initial data for the einstein equations, 2021.
\newblock arXiv:2002.01534.

\bibitem{Hirsch-Kazaras-Khuri-Zhang}
Sven Hirsch, Demetre Kazaras, Marcus Khuri, and Yiyue Zhang.
\newblock Rigid comparison geometry for {R}iemannian bands and open incomplete
  manifolds, 2022.
\newblock arXiv:2209.12857.

\bibitem{HZ}
Sven Hirsch and Yiyue Zhang.
\newblock Stability of llarull's theorem in all dimensions, 2023.
\newblock arXiv:2310.14412.

\bibitem{Huisken-Ilmanen}
Gerhard Huisken and Tom Ilmanen.
\newblock The inverse mean curvature flow and the {R}iemannian {P}enrose
  inequality.
\newblock {\em J. Differential Geom.}, 59(3):353--437, 2001.

\bibitem{Sweeney}
Paul~Sweeney Jr.
\newblock Examples for scalar sphere stability, 2023.
\newblock arXiv:2301.01292.

\bibitem{Lee-Kazaras-Khuri-PMT-Stability}
Demetre Kazaras, Marcus Khuri, and Dan Lee.
\newblock Stability of the positive mass theorem under ricci curvature lower
  bounds, 2021.
\newblock arXiv:2111.05202.

\bibitem{KK}
Demetre Kazaras and Kai Xu.
\newblock Drawstrings and flexibility in the georch conjecture, 2023.
\newblock arXiv:2309.03756.

\bibitem{Lee-PMT-Conformal}
Dan~A. Lee.
\newblock On the near-equality case of the positive mass theorem.
\newblock {\em Duke Math. J.}, 148(1):63--80, 2009.

\bibitem{LeeSormani1}
Dan~A. Lee and Christina Sormani.
\newblock {S}tability of the positive mass theorem for rotationally symmetric
  riemannian manifolds.
\newblock {\em Journal fur die Riene und Angewandte Mathematik (Crelle's
  Journal)}, 686, 2014.

\bibitem{LNN}
Man-Chun Lee, Aaron Naber, and Robin Neumayer.
\newblock {$d_p$}-convergence and {$\epsilon$}-regularity theorems for entropy
  and scalar curvature lower bounds.
\newblock {\em Geom. Topol.}, 27(1):227--350, 2023.

\bibitem{Lee-Topping}
Man-Chun Lee and Peter~M. Topping.
\newblock Metric limits of manifolds with positive scalar curvature, 2022.
\newblock arXiv:203.01223.

\bibitem{Llarull}
Marcelo Llarull.
\newblock Sharp estimates and the {D}irac operator.
\newblock {\em Mathematische Annalen}, pages 55--71, 1998.

\bibitem{PKP19}
Armando J.~Cabrera Pacheco, Christian Ketterer, and Raquel Perales.
\newblock Stability of graphical tori with almost nonnegative scalar curvature.
\newblock {\em Calc. Var.}, 59(134), 2020.

\bibitem{Schoen-Yau-min-surf}
R.~Schoen and Shing~Tung Yau.
\newblock Existence of incompressible minimal surfaces and the topology of
  three-dimensional manifolds with nonnegative scalar curvature.
\newblock {\em Ann. of Math. (2)}, 110(1):127--142, 1979.

\bibitem{SchoenYauPMT}
Richard Schoen and Shing~Tung Yau.
\newblock On the proof of the positive mass conjecture in general relativity.
\newblock {\em Comm. Math. Phys.}, 65(1):45--76, 1979.

\bibitem{Stavrov-Sormani}
Christina Sormani and Iva Stavrov~Allen.
\newblock Geometrostatic manifolds of small {ADM} mass.
\newblock {\em Comm. Pure Appl. Math.}, 72(6):1243--1287, 2019.

\bibitem{SW-JDG}
Christina Sormani and Stefan Wenger.
\newblock The intrinsic flat distance between {R}iemannian manifolds and other
  integral current spaces.
\newblock {\em J. Differential Geom.}, 87(1):117--199, 2011.

\bibitem{Witten-PMT}
Edward Witten.
\newblock A new proof of the positive energy theorem.
\newblock {\em Comm. Math. Phys.}, 80(3):381--402, 1981.

\end{thebibliography}






\end{talk}

\end{document}